 \documentclass[11pt]{article}
 \usepackage[top=1in,bottom=1in,left=1.5in,right=1.5in]{geometry}
  \usepackage{amsmath,amssymb}
  \usepackage{latexsym}
  \usepackage[dvips]{pstricks} 
   \usepackage{pst-node}
  \usepackage[dvips]{graphicx}
  
  \newcommand{\C}{\mathbb{C}}
  \newcommand{\F}{\mathbb{F}}
  \newcommand{\N}{\mathbb{N}}
  
  \newcommand{\Q}{\mathbb{Q}}
  \newcommand{\R}{\mathbb{R}}

  \newcommand{\e}{\mathbf{e}}

  \newcommand{\U}{\mathbf{U}}
  \newcommand{\uu}{\mathbf{u}}
  \newcommand{\vv}{\mathbf{v}}
  \newcommand{\V}{\mathbf{V}}
  \newcommand{\w}{\mathbf{w}}
  
  \newcommand{\x}{\mathbf{x}}
  
  \newcommand{\y}{\mathbf{y}}
  
  \newcommand{\z}{\mathbf{z}}
  \newcommand{\0}{\mathbf{0}}

  \newcommand{\cA}{\mathcal{A}}

  \newcommand{\cS}{\mathcal{S}}
  \newcommand{\cT}{\mathcal{T}}

  \newcommand{\rO}{\mathrm{O}}

  \newcommand{\lan}{\langle}
  \newcommand{\ran}{\rangle}
  \newcommand{\an}[1]{\lan#1\ran}
  \def\diag{\mathop{{\rm diag}}\nolimits}
  \newcommand{\hs}{\hspace*{\parindent}}
  \newcommand{\proof}{\hs \textbf{Proof.\ }}
  
  \newcommand{\tr}{\mathop{\mathrm{tr}}\nolimits}

  \newcommand{\Sym}{\mathop{\mathrm{Sym}}\nolimits}
  \newcommand{\trans}{^\top}
  \newcommand{\qed}{\hspace*{\fill} $\Box$\\}

  \newcommand{\rS}{\mathrm{S}}

  \newcommand{\rank}{\mathrm{rank\;}}

  \newtheorem{theo}{\bfseries \hs Theorem}
  
  \newtheorem{prop}[theo]{\bfseries \hs Proposition}
  
  \newtheorem{lemma}[theo]{\bfseries \hs Lemma}
  \newtheorem{corol}[theo]{\bfseries \hs Corollary}
  \newtheorem{con}[theo]{\bfseries \hs Conjecture}

  \numberwithin{equation}{section} 

 \renewcommand{\span}{\mathrm{span}}

\begin{document}

 \title{Best rank one approximation of real symmetric tensors\\
 can be chosen symmetric\footnotemark[1]}

 \author{
 S. Friedland\footnotemark[2]\\
 Department of Mathematics, Statistics and Computer Science\\
 University of Illinois at Chicago\\
 Chicago, Illinois 60607-7045, USA\\
 e-mail:\texttt{friedlan@uic.edu}}

 \renewcommand{\thefootnote}{\arabic{footnote}}
 \footnotetext[1]{ To appear in Special Issue ``Tensor Theory'', \textit{Frontiers of Mathematics in China}, Springer.}
 \footnotetext[2]{This work was supported by NSF grant DMS-1216393.}
 \date{November 25, 2012 }
 \maketitle
 \begin{abstract}
 We show that a best rank one approximation to a real symmetric tensor, which in principle can be nonsymmetric, can be chosen symmetric.
 Furthermore, a symmetric best rank one approximation to a symmetric tensor is unique if the tensor does not lie on a certain
 real algebraic variety.

 \end{abstract}



 \renewcommand{\thefootnote}{\arabic{footnote}}

 \section{Introduction}\label{intro}
 Denote by $\R^{n_1\times\ldots\times n_d}:=\otimes_{i=1}^d \R^{n_j}$
 the tensor products of $\R^{n_1},\ldots,\R^{n_d}$ for an integer $d\ge 2$.
 $\cT=[t_{i_1,\ldots, i_d}]\in \R^{n_1\times\ldots\times n_d}$
 is called a \emph{$d$-tensor}.  Let $[d]:=\{1,\ldots,d\}$.
 For $\x_j=(x_{1,j},\ldots,x_{n_j,j})\trans\in\R^{n_j}, j\in [d]$, the decomposable tensor $\otimes_{j=1}^d \x_j=\x_1\otimes\ldots\x_d$
 is given as $[\prod_{j=1}^d x_{i_j,j}]\in\R^{n_1\times\ldots\times n_d}$.  A decomposable tensor is a rank one tensor if and only if $\x_j\ne \0$ for each $j\in[d]$.

 On $\R^{n_1\times\ldots\times n_d}$ define the inner product $\an{\cS,\cT}:=\sum_{i_j\in[n_j],j\in[d]} s_{i_1,\ldots,i_d}t_{i_1,\ldots,i_d}$
 and the Hilbert-Schmidt norm $\|\cT\|:=\sqrt{\an{\cT,\cT}}$.  For $\x\in\R^n$ let $\|\x\|:=\sqrt{\x\trans \x}$ be the $\ell_2$ norm.
 Denote by $\rS^{n-1}:=\{\x\in\R^n,\;\|\x\|=1\}$ the $n-1$ dimensional unit sphere.
 Then a best rank one approximation of $\cT$ is a decomposable tensor solving the minimal problem
 \begin{equation}\label{brank1appr}
 \min_{s\in\R,\x_j\in\rS^{n_j-1},j\in[d]} \|\cT-s\otimes_{j=1}^d\x_j\|=\|\cT-a\otimes_{j=1}^d\uu_j\|.
 \end{equation}
 It is well known
 \begin{equation}\label{maxmultform}
 a=\max_{\x_j\in\rS^{n_j-1},j\in[d]}\sum_{i_j\in[n_j],j\in[d]}t_{i_1,\ldots,i_d} x_{i_1,1}\ldots x_{i_d,d}=\an{\cT,\otimes_{j=1}^d \uu_j}.
 \end{equation}

 For matrices, i.e. $d=2$, a best rank one approximation of $A\in\R^{m\times n}$ is given by $\sigma_1(A)\x\y\trans, \x\x\trans=\y\y\trans=1$,  where $\sigma_1(A)$ is the maximal singular value of $A$, and $A\y=\sigma_1(A)\x, A\trans \x=\sigma_1(A)\y$.
 A best rank one approximation is unique if and only if $\sigma_1(A)>\sigma_2(A)$.  Assume that $A\in\R^{n\times n}$ is symmetric.
 Let $\lambda_1(A)\ge\ldots\ge\lambda_n(A)$ be the $n$-eigenvalues of $A$, counted with their multiplicities.
  Then $\sigma_1(A)=\max(|\lambda_1(A), |\lambda_n(A)|)$.  Furthermore, there always exists a symmetric best rank approximation.  Moreover, there exists a nonsymmetric best rank one approximation if and only if $\lambda_n(A)=-\lambda_1(A)<0$.

 Assume now that $\cT$ is a $d$-mode tensor with $d\ge 3$.
 Let $\alpha$ be a subset of $[d]$ of cardinality $2$ at least.
 We say that $\cT=[t_{i_1,\ldots,i_d}]\in \R^{n_1\times\ldots n_d}$ is symmetric with respect to $\alpha$ if $n_p=n_q$ for each pair $\{p,q\}\subset \alpha$  and the value of $t_{i_1,\ldots,i_d}$ does not change if we interchange any two indices $i_p,i_q$ for $p,q\in\alpha$ and for any choice $i_j\in [n_j], j\in [d]$.  We agree that any tensor $\cT$ is symmetric with respect to each $\{i\}$, for $i\in [d]$.
 $\cT$ is called symmetric if $\alpha=[d]$.
 Denote by $\Sym(n,d)$ all $d$-mode symmetric tensors in $\R^{n\times\ldots\times n}$.  These tensors are called also supersymmetric.
 Let $\cT\in \R^{n_1\times \ldots\times n_d}$ be given.
 Clearly, there exists a unique decomposition of $[d]$ to a disjoint union of nonempty sets $\cup_{j=1}^m \alpha_j$ such that the following
 conditions hold.
 \begin{itemize}
 \item For each $j\in [m]$ the tensor $\cT$ is symmetric with respect to $\alpha_j$.
 \item For $1\le j <k\le m$ and two indices $p\in\alpha_j,q\in\alpha_k$ the tensor $\cT$ is not symmetric with respect to $\{p,q\}$.
 \end{itemize}
 We call $[d]=\cup_{j=1}^m\alpha_j$ the \emph{symmetric decomposition} for $\cT$.

 The main result of this paper is.
 \begin{theo}\label{maintheo}  Let $\cT\in  \R^{n_1\times\ldots n_d}\setminus\{0\}$ be given.  Assume that $[d]=\cup_{j=1}^m \alpha_j$ is the symmetric
 decomposition for $\cT$.  Then there exist a best rank one approximation $\cA$ to $\cT$ such that $\cA$
 is symmetric with to each $\alpha_j$.
 \end{theo}

 In the special case of symmetric tensors, i.e. $m=1$, the above theorem is also proved in \cite[Theorem 4.1]{CHLZ12}.

 It is not difficult to show that to prove the above theorem it is enough to show the above theorem for symmetric tensors.
 Furthermore, it is enough to show the above theorem for $\cT\in\Sym(2,d)$.  Finally we show that there exists an real algebraic variety
 $\Sigma_1(n,d)\subset \Sym(n,d)$ such that for $\cT\in \Sym(n,d)\setminus\Sigma_1(n,d)$ a best rank one symmetric approximation is unique.

 We now describe briefly the contents of our paper.  In \S2 we summarize the well known results of best rank one approximation for real
 matrices that are used in this paper.  In \S3 we discuss certain basic results on best rank one approximation of real tensors.
 In \S4 we give a complete characterization of tensors $\cT\in\Sym(2,3)$ which have nonsymmetric best rank one approximation.
 In \S5 we prove Theorem \ref{maintheo} by first showing the case where $\cT$ is a symmetric tensor, (Theorem \ref{brank1symtheo}).
 In \S6 we show that a ''generic" symmetric tensor has a unique best rank one approximation.  The last section uses some facts from
 algebraic geometry, and probably the most difficult section of this paper.

 \section{Best rank one approximation of matrices}

 We recall briefly the needed results on best rank one approximation of matrices $A=[a_{i,j}]\in\R^{m\times n}$ \cite{GolV96}.
 Assume that $r=\rank A$.  Then $A$ has exactly $r$ positive singular values $\sigma_1(A)\ge\ldots\ge \sigma_r(A)>0$.
 $\sigma_1(A)^2\ge\ldots\ge\sigma_r(A)^2$ are the positive eigenvalues of either $AA\trans$ or $A\trans A$.
 \begin{equation}\label{charsingval1}
 \sigma_1(A)=\max_{\x\in\rS^{m-1},\y\in\rS^{n-1}}\x\trans A\y.
 \end{equation}
 The left and the right singular pair of singular vectors $\uu\in\rS^{m-1},\vv\in\rS^{n-1}$ corresponding to $\sigma_1(A)$ are given by equalities
 \begin{equation}\label{uvsingvectA}
 A\trans\uu=\sigma_1(A)\vv,\; A\vv=\sigma_1(A)\uu, \quad \uu\in\rS^{m-1},\vv\in\rS^{n-1}.
 \end{equation}
 Hence $\sigma_1(A)=\uu\trans A\vv$.  Furthermore any best rank one approximation of $A$ in the Frobenius norm $\|A\|=\sqrt{\tr (A\trans A)}$
 is of the form $\sigma_1(A)\uu\trans \vv$ for some pair of singular vectors $\uu,\vv$ corresponding to $\sigma_1(A)$:
 \begin{equation}\label{brank1apmat}
 \min_{s\in \R,\x\in\rS^{m-1},\y\in\rS^{n-1}}\|A-s\x\trans\y\|=\|A-\sigma_1(A)\uu\trans\vv\|.
 \end{equation}

 Recall that $\Sym(n,2)\subset \R^{n\times n}$ is the space of symmetric matrices.   The following result is well known and we bring its proof for completeness.
 \begin{prop}\label{branapproxsymmat}  Let $A\in\Sym(n,2)$.  Then $\sigma_1(A)=\max(|\lambda_1(A)|,|\lambda_n(A)|)$ and $A$ has a symmetric best rank one approximation. Suppose furthermore that $A\ne 0$.
 Then
 \begin{enumerate}
 \item\label{symrank1motappr} Any best rank one approximation to $A$ is symmetric if and only if $\lambda_1(A)\ne -\lambda_n(A)$.
 \item\label{nsmrank1motappr} Assume that $\lambda_1(A)=-\lambda_n(A)$.  Then $\sigma_1(A)=\lambda_1(A)$.  Furthermore, $\sigma_1(A)\uu\trans\vv$, where \eqref{uvsingvectA} holds,
 is a nonsymmetric best rank one approximation of $A$ if neither $\uu$ nor $\vv$ are eigenvectors of $A$ and $\uu,\vv$ are eigenvectors of
 $A^2$ corresponding to $\lambda_1(A)^2$.
 \end{enumerate}
 \end{prop}
 \proof  Since $A\trans A=A^2$ the singular values of $A$ are $|\lambda_i(A)|, i\in [n]$.
 As all eigenvalues of $A$ lie in the interval $[\lambda_n(A),\lambda_1(A)]$ it follows that
 $\sigma_1(A)=\max_{i\in[n]}|\lambda_i(A)|=\max(|\lambda_1(A)|,|\lambda_n(A)|)$.
 We now show that there exists a symmetric best rank one approximation of $A$.  Clearly, it is enough to consider $A\ne 0$.
 Assume that $\sigma_1(A)=|\lambda_j(A)|, j\in \{1,n\}$.  Let $\vv\in \rS^{n-1}$ be an eigenvector of $A$ corresponding to $\lambda_j(A)$.
 Hence $\uu=\frac{\lambda_j(A)}{|\lambda_j(A)|}\vv$ and $\sigma_1(A)\uu\trans\vv=\lambda_j(A)\vv\trans\vv$, which is a symmetric best rank one
 approximation of $A$.

 Assume now that $A\ne 0$.  Note that the assumption that $|\lambda_1(A)|\ne |\lambda_n(A)|$ is equivalent to $\lambda_1(A)\ne -\lambda_n(A)$.
 Assume first that $\lambda_1(A)\ne -\lambda_n(A)$.  Then there exists a unique $j\in \{1,n\}$ such that $\sigma_1(A)=|\lambda_j(A)|$.
 As a right singular vector $\vv\in\rS^{n-1}$ of $A$ is an eigenvector of $A^2$ corresponding to $\lambda_j(A)^2$ it follows that
 $\vv$ is an eigenvector $A$.  Hence the above argument show that $\sigma_1(A)\uu\trans \vv$ is symmetric and \ref{symrank1motappr} holds.

 Assume now that $\lambda_1(A)=-\lambda_n(A)$.  If $\vv$ is an eigenvector of $A$ corresponding to $\lambda_j(A), j\in\{1,n\}$ then $\uu$
 is also eigenvector of $A$ corresponding to $\lambda_j(A)$, and vice versa.  In this case $\sigma_1(A)\uu\trans \vv$ is a symmetric best rank one approximation.
 Assume that $\vv$ is not an eigenvector of $A$.  Since $\vv$
 is a right singular vector corresponding to $\sigma_1(A)$ it follows that $\vv$ is an eigenvector of $A^2$ corresponding to the eigenvalue
 $\lambda_j(A)^2$.  Similar claim holds for $\uu$.  Hence $(A\vv)\trans \vv$ is a nonsymmetric best rank approximation
 and \ref{nsmrank1motappr} holds.  \qed

  \begin{lemma}\label{uniqaprmat}

 \noindent
 \begin{enumerate}
 \item \label{uniqaprmatnsym} $A\in\R^{m\times n}$ has a unique rank one approximation for $m,n\ge 2$, unless $A$ lies on a subvariety of codimension two.
 \item \label{uniqaprmatsym}
 $A\in\Sym(n,2)$ has a unique rank one approximation for $n\ge 2$, which is symmetric, unless $A$ lies on a subvariety of codimension one in $\Sym(n,2)$.
 \end{enumerate}
 \end{lemma}
 \proof
 To prove the first part of the Lemma we use the singular value decomposition.  Assume without loss of generality that $2\le m \le n$.
 Then each matrix $A\in\R^{m\times n}$ is of the form $UDV\trans$ where $U=[\uu_1,\ldots\uu_m]\in \R^{m\times m}$ is orthogonal,
 $V=[\vv_1,\ldots,\vv_m]\in \R^{n\times m}$ has $m$-orthogonal columns and $D=\diag(d_1,\ldots,d_m)$,, where $d_1\ge\ldots\ge d_m\ge 0$.  So $\sigma_i(A)=d_i, i\in[m]$, and the columns $i$
 of $V$ and $U$ are right and the left singular vectors respectively corresponding to the singular value $d_i$.
 Note that if $d_1>\ldots>d_m>0$  then each column of $V$ is determined up to $\pm$, and after $V$ is fixed, then $U$ is determined uniquely.
 In this case $\sigma_1(A)=d_1>\sigma_2(A)=d_2$ and $A$ has a unique rank one approximation.
 $A$ has no unique approximation if and only if $d_1=d_2$.  The generic case for this situation is
 \begin{equation}\label{genmatdegcase}
 d_1=d_2>d_3>\ldots>d_m>0.
 \end{equation}
 The equality $d_1=d_2$ means that we loose one parameter.  The columns $3,\ldots,m$ of $V$ are determined uniquely to $\pm 1$.
 The first two columns of $V$ are not determined uniquely.  What is determined uniquely is the two dimensional subspace $\V\subset\R^n$
 which is orthogonal the the columns $3,\ldots,m$ of $V$ and the null space of $A$.  We can choose as a first column $\vv_1$ of $V$ any unit vector
 in $\V$.  The the second column $\vv_2$ of $V$ is a unit vector in $\V$ which is orthogonal to $\vv_1$.  So $\vv_2$ is determined uniquely up to a sign.  Recall that  $\uu_i=d_i^{-1}A\vv_i,i\in[n]$.  Hence the set $A=UDV\trans$ of the above form, where the entries of $D$ satisfy \eqref{genmatdegcase}, is a manifold $\Phi(m,n)\subset \R^{m\times n}$ of codimension two in $\R^{m\times n}$.  It is left to show that there is a nonzero polynomial $Q$ in the entries of $A$ which satisfies the following conditions.  First, $Q$ vanishes on $\Phi(m,n)$.  Second, for each  $A\in \Phi(m,n)$ the exists a neighborhood $\rO\subset \R^{m\times n}$ of $A$ such that the zero set of $Q$
 on $\rO$ is equal to $\rO\cap\Phi(m,n)$.  Consider the symmetric matrix $B=AA\trans$.  Let $\textrm{Dis}(B)$ be the discriminant of the characteristic polynomial of $B$.  Then $\textrm{Dis}(B)$ vanishes if and only if $B$ has a multiple eigenvalue.  In particular $\textrm{Dis}(B)$ vanishes on $\Phi(m,n)$.  Fix $A\in \Phi(m,n)$.  Assume that $C\in \R^{m\times n}$ is very close to $A$.  Then $\sigma_2(C)>\ldots>\sigma_m(C)>0$.
 So $Q(CC\trans)=0$ if and only if $\sigma_1(C)=\sigma_2(C)$.  This establishes the first part of the Lemma.

 Recall that the set of all $A\in\Sym(n,2)$ having at least one multiple eigenvalue is a variety $\Delta_n$ of codimension two, e.g. \cite{FRS}.
 (This follows from the arguments of the first part of the Lemma.  $\Delta_n$ is the zero set of the discriminant of the characteristic polynomial of $A$.)   Assume that $A\in \Sym(n,2)\setminus\Delta_n$. (So $A\ne 0$.)  Proposition \ref{branapproxsymmat} yields that $A$ has a unique symmetric rank one approximation if and only if $\lambda_1(A)\ne -\lambda_n(A)$.  It is left to show that all matrices $A\in\Sym(n,2)\setminus\Delta_n$ satisfying $\lambda_1(A)+\lambda_n(A)=0$ lie on a variety of codimension one.
 This follows from the spectral decomposition $A=UDU\trans$, where $U$ ia an orthogonal matrix and $D=\diag(d_1,\ldots,d_n),d_1>\ldots>d_n$.
 Note that the columns of $U\trans$ are determined uniquely up to a sign.
 Hence the set of all $A\in\Sym(n,2)\setminus\Delta_n$ satisfying $\lambda_1(A)+\lambda_n(A)=0$ is a manifold of codimension one.  It is left to show that that
 this manifold is a zero set of some polynomial in $A$.  Recall that for $A\in\Sym(n,2)$ the matrix $A\otimes I_n+I_n\otimes A\in\Sym(n^2,2)$,
 where $I_n\in\Sym(n,2)$ is the identity matrix and $A\otimes B$ is the Kronecker tensor product, have eigenvalues $\lambda_i(A)+\lambda_j(A)$ for $i,j\in[n]$.  Hence the zero set of $\det(A\otimes I_n+I_n\otimes A)$ includes the above manifold.
 This concludes the proof of the second part of the Lemma.
 \qed

 \section{Preliminary results on best rank one approximation of tensors}

 Recall that $\infty$-Schatten norm of $\cT=[t_{i_1,\ldots,i_d}]\in \R^{n_1\times \ldots\times n_d}$, with respect to the $\ell_2$ norm on each factor $\R^{n_i}$, is given by
 \[\|\cT\|_{\infty,2}:=\max_{\x_i\in\rS^{n_i-1}, i\in [d]}|\an{\cT\times \otimes_{i=1}^d \x_i}|.\]
 Since $-\rS^{n-1}=\rS^{n-1}$ it follows that
 \begin{equation}\label{definfsnorm}
 \|\cT\|_{\infty,2}:=\max_{\x_i=(x_{1,i},\ldots,x_{n_i,i})\trans\in\rS^{n_i-1}, i\in [d]} \sum_{i_1=\ldots=i_d=1}^{n_1,\ldots,n_d}t_{i_1,\ldots,i_d}x_{i_1,1}\ldots x_{i_d,d}.
 \end{equation}
 See for example \cite{DF93} for a modern exposition on tensor norms and \cite{Fr82} for simple geometrical properties of cross norms.
 Note that for matrices, i.e. $d=2$, $\|A\|_{\infty,2}$ is the operator norm $\|A\|_2=\sigma_1(A)$, where $A\in\R^{m\times n}$ viewed as
 a linear operator $\y\mapsto A\y$ from $\R^n$ to $\R^m$.

 Let $\beta\subset[d]$ be a nonempty set and assume that $\x_j=(x_{1,j},\ldots,x_{n_j,j})\trans\in \R^{n_j}$ for $j\in\beta$.
 Denote by $\cT\times \otimes_{j\in\beta}\x_j$ the contracted $d-|\beta|$ tensor
 \begin{equation}\label{defTcontr}
 \cT\times\otimes_{j\in\beta}\x_j:=\sum_{i_j\in[n_j],j\in\beta}t_{i_1,\ldots,i_d}\prod_{j\in\beta} x_{i_j,j} \in \otimes_{k\in[d]\setminus\beta}\R^{n_k}.
 \end{equation}
 Note that if $\beta=[d]$ then $\cT\times\otimes_{j\in[d]}\x_j=\an{\cT,\otimes_{j\in[d]}\x_j}$.
 Let $\beta$ be a nonempty strict subset of $[d]$.   By considering the maximum in \eqref{defTcontr} as a maximum on $\x_j$ first on $j\in[d]\setminus\beta$ and then on $j\in\beta$ we deduce
 \begin{equation}\label{infshatnrmchar}
 \|\cT\|_{\infty,2}=\max_{ \x_j\in \rS^{n_j-1},j\in\beta}\|\cT\times\otimes_{j\in\beta}\x_j\|_{\infty,2}.
 \end{equation}
 Suppose that $\beta=[d]\setminus\{p,q\}$, where $1\le p <q\le d$.  We view $\cT\times \otimes_{j\in\beta}\x_j$ as a matrix in $\R^{n_p\times n_q}$. Hence
 \begin{equation}\label{infshatnrmchar1}
 \|\cT\|_{\infty,2}=\max_{ \x_j\in \rS^{n_j-1},j\in[d]\setminus\{p,q\}}\sigma_1(\cT\times\otimes_{j\in[d]\setminus\{p,q\}}\x_j).
 \end{equation}

 The following result is well known and we bring its proof for completeness.
 \begin{lemma}\label{brank1char}  Let $\cT\ne 0$ be a given tensor in $\R^{n_1\times \ldots\times n_d}$.
 Then $a\otimes_{i=1}^d \uu_i$, where $\uu_i\in\rS^{n_i-1}, i\in[d]$, is a best rank one approximation of $\cT$ if and only if the following conditions hold.  First $a=\pm\|\cT\|_{\infty,2}$.   Second the function
 $\an{\cT,\otimes_{j\in[d]}\x_j}$ attains its maximum or minimum on $\rS^{n_1-1}\times\ldots\times\rS^{n_d-1}$ at $(\uu_1,\ldots,\uu_d)$.
 In particular
 \begin{equation}\label{brank1char1}
 \cT \times \otimes_{j\in[d]\setminus\{i\}} \uu_j=\lambda \uu_i, \quad i\in [d],
 \end{equation}
 where $\lambda=\pm\|\cT\|_{2,\infty}$ and $\uu_i\in\rS^{n_i-1}$ for $i\in[d]$.
 Suppose furthermore that $\cT$ is symmetric with respect to $1\le p <q\le d$.  Then there exist a best rank one approximation
 which is symmetric with respect with respect to $\{p,q\}$.
 \end{lemma}
 \proof
  Let $\x_i\in\rS^{n_i-1}$ for $i\in[d]$.  Then $\otimes_{i\in[d]}\x_i$ is a unit vector
 in $\otimes_{i=1}^d \R^{n_i}$.  Let $\U:=\span(\otimes_{i\in[d]}\x_i)$ and $\U^{\perp}$ be the orthonormal complement of $\U$ in $\otimes_{i=1}^d \R^{n_i}$.  The orthogonal projection of $\cT$ on $\U$ is given by
 $P_{\U}(\cT)=\an{\cT,\otimes_{i\in[d]}\x_i}\otimes_{i\in[d]}\x_i$, and $\|P_{\U}(\cT)\|=|\an{\cT,\otimes_{i\in[d]}\x_i}|$.  It is well known that $\min_{s\in \R} \|\cT-s\otimes_{i=1}^d \x_i\|=
 \|P_{\U^{\perp}}(\cT)\|$.
 The Pythagoras theorem yields that
 \begin{equation}\label{pythid}
 \|\cT\| = \|P_{\U}(\cT)\|^2 + \|P_{\U^\perp}(\cT)\|^2=\an{\cT,\otimes_{i\in[d]}\x_i}^2+\|P_{\U^\perp}(\cT)\|^2.
 \end{equation}
 Hence a minimal solution of the left-hand side of \eqref{brank1appr} gives rise to a maximum or minimum of
 $\an{\cT,\otimes_{j\in[d]}\x_j}$ on $\rS^{n_1-1}\times\ldots\times\rS^{n_d-1}$.

 We now give a short proof of a result by Lim \cite{Lim05}.
 Consider the maximum problem \eqref{maxmultform}.  Use Lagrange multipliers for the function   $\an{\cT,\otimes_{j\in[d]}\x_j}-\sum_{j\in[d]}\lambda_j\x_j\trans\x_j$ to deduce that a maximum solution satisfies
 \[\cT\times\otimes_{j\in[d]\setminus\{i\}}\uu_j=\lambda_i\uu_i, \quad i\in[d].\]
 Hence $\an{\cT,\otimes_{j\in[d]}\uu_j}=\lambda_i\uu_i\trans\uu_i=\lambda_i$ for each $i\in[d]$.  Therefore a best rank one approximation
 $a\otimes_{j\in[d]}\uu_j$ satisfies $a=\|\cT\|_{\infty,2}$ and \eqref{brank1char1}, where $\lambda=a$ and $\uu_j\in\rS^{n_j-1}$ for $j\in[d]$.
 Similar results hold for the minimum of $\an{\cT,\otimes_{j\in[d]}\x_j}$ on $\rS^{n_1-1}\times\ldots\times\rS^{n_d-1}$.

 Assume now that $\cT$ is symmetric with respect to two indices $p<q$. So $n_p=n_q$.  Assume that $a\otimes_{i=1}^d \uu_i$ is best rank one approximation.  Let $\U_i=\span(\uu_i), i\in [d]$.  Note that $A:=\cT\times \otimes_{i\in [d]\setminus\{p,q\} }\uu_i$ is a symmetric matrix.
 As best rank one approximation of a $A$ can be chosen symmetric we deduce that we can choose $\uu_p,\uu_q\in \rS^{n_p-1}$ such that
 $\uu_q\in \{\uu_p,-\uu_p\}$.  Hence there exist a best rank approximation of $\cT$ which is symmetric with respect to $\{p,q\}$.
 \qed

 \begin{corol}\label{charb1bap}  Let $\cT\ne 0$ be a given tensor in $\R^{n_1\times \ldots\times n_d}$.
 Then rank one tensor $\cA\in\R^{n_1\times \ldots\times n_d}$ is best rank one approximation of $\cT$ if and only if
 \begin{equation}\label{charb1bap1}
 \an{\cT,\cA}=\|\cT\|_{\infty,2}^2=\|\cA\|^2.
 \end{equation}
 \end{corol}

 The following lemma is straightforward.
 \begin{lemma}\label{uniqbrankap}  Let $\cT\in\R^{n_1\times\ldots\times n_d}$ and assume that $a\otimes_{j\in[d]}\uu_j$ is a best rank one approximation of $\cT$.
 Suppose that $\cT$ is symmetric with respect to $\alpha\subset[d]$.
 Let $\sigma:[d]\to[d]$ be a permutation which is identity on $[d]\setminus\alpha$.  Then $a\otimes_{j\in[d]} \x_{\sigma(j)}$ is a best
 is rank one approximation of $\cT$.  In particular, if $a\otimes_{j\in[d]}\uu_j$ is unique best rank one approximation of $\cT$
 then $a\otimes_{j\in[d]}\uu_j$ is symmetric with respect to $\alpha\subset[d]$.
 \end{lemma}

 Lemma \ref{uniqaprmat} suggests the following conjecture.
 \begin{con}\label{brank1con} Let $d\ge 3, n_j\ge 2, j\in [d]$ be integers.  Then
 \noindent
 \begin{enumerate}
 \item \label{uniqaprtengen} $\cT\in\R^{n_1\times\ldots\times n_d}$ has a unique rank one approximation, unless $\cT$ lies on a subvariety.
 \item \label{uniqaptensym}
 $\cT\in\Sym(n,d)$ has a unique rank one approximation, which is symmetric, unless $\cT$ lies on a subvariety.
 \end{enumerate}

 \end{con}

 The above conjecture was recently settled in \cite{FO12}. 
 We remark that T. Kolda in her lecture \cite{Kol} stated a stronger version of the second part of Conjecture \ref{brank1con}.
 Namely:''rank-r symmetric factorization of a symmetric tensor is unique even without the symmetry condition, under very mild conditions",
 although she did not specify the nature of the mild conditions.

 \begin{prop}\label{conjexistsymbr1ap}  Assume that the second part of Conjecture \ref{uniqaptensym}
 holds for some integers $n\ge 2, d\ge 3$.  Then each $\cT\in\Sym(n,d)$ has a best rank one symmetric approximation.
 \end{prop}
 \proof  Assume that $\Psi_n\subset \Sym(n,d)$ is the variety of all symmetric tensors which do not have a unique best rank one approximation.
 So each $\cT\in\Sym(n,d)\setminus\Psi_n$ has a unique best rank one approximation $\cA(\cT)$ which is symmetric.  Assume now that $\cT\in\Psi_n$.
 As $\Psi_n$ is a variety, there exists a sequence $\cT_k, k\in\N$ which converges to $\cT$, and $\cT_k\not\in\Psi_n$ for all $k\in\N$.
 Use Corollary \ref{charb1bap} to deduce
 \[ \|\cT\|_{\infty,2}^2=\lim_{k\to\infty}\|\cT_k\|_{2,\infty}^2=\lim_{k\to\infty}\an{\cT_k,\cA(\cT_k)}=\lim_{k\to\infty} \|\cA(\cT_k)\|^2.\]
 So $\cA(\cT_k)$ is a bounded sequence in $\Sym(n,d)$.  Hence there exists a subsequence $\cA(\cT_{k_l})$ which converges to a rank
 one symmetric tensor $\cA$ which satisfies $\|\cT\|_{\infty,2}^2=\an{\cT,\cA}=\|\cA\|^2$.  Hence by Corollary \ref{charb1bap} $\cA$ is a symmetric best rank one approximation of $\cT$.  \qed

 A weaker version of the second part of Conjecture \ref{brank1con} is:

 \begin{theo}\label{brank1symtheo}  Every symmetric tensor $\cT\in\Sym(n,d)$ has a symmetric best rank one approximation for integers $n\ge 2, d\ge 3$.
 \end{theo}
 Note that the above theorem is a special case of Theorem \ref{maintheo}.  We will first prove Theorem \ref{brank1symtheo}, and using it we will prove Theorem
 \ref{maintheo}.

 \begin{lemma}\label{nto2}  Let $k\ge 2$.  Assume that Theorem \ref{brank1symtheo} holds for $n=2$ and for all positive integers $d$ in the interval $[2,k]$.  Then Theorem \ref{brank1symtheo}
 holds for all integers $n\ge 3$ and $d\in [2,k]$.
 \end{lemma}
 \proof  We prove our Lemma by induction on $k$.
 In view of Proposition \ref{branapproxsymmat} Theorem \ref{brank1symtheo} trivially holds for $k=2$.  Assume that $N\ge 3$ and suppose that we proved the Lemma for $k=N-1$. Assume that Theorem \ref{brank1symtheo} holds for $\Sym(2,N)$.  Let $\cT=[t_{i_1,\ldots,i_N}]\in \Sym(n,N)$ and $n\ge 3$.  Suppose that $\|\cT\|_{2,\infty}=|\an{\cT,\otimes_{j\in[d]}\vv_j}|$, where $\vv_j\in \rS^{n-1}$ for $j\in [N]$.  Let $\cS:=\cT\times \vv_N\in \Sym(n,N-1)$.
 So $\|\cS\|_{2,\infty}=\|\cT\|_{2,\infty}$.
 Our induction assumption implies that there exists $\uu\in\rS^{n-1}$ such that $\|\cS\|_{2,\infty}=|\an{\cS,\otimes_{j\in[N-1]}\uu_j}|$,
 where $\uu_j=\uu$ for $j\in[N-1]$.  Let $\uu_N=\vv_N$.  Then $\|\cT\|_{\infty,2}=|\an{\cT,\otimes_{j\in[N]}\uu_j}|$.
 If $\uu_N=\pm \uu$ it follows that that $\an{\cT,\otimes_{j\in[N]}\uu_j}\otimes_{j\in[N]}\uu_j$ a best rank one symmetric approximation
 of $\cT$, and we are done.  Suppose that $\uu_N\ne \pm \uu$.  So $\span(\uu,\uu_N)$ is two dimensional.
 By changing an orthonormal basis in $\R^n$ we may assume without loss of generality that $\span(\uu,\uu_N)=\span(\e_1,\e_2)$, where
 $\e_j$ is the $j$-th column of the identity matrix $I_n$.  Let $\cT'=[t_{i_1,\ldots,i_N}]_{i_1,\ldots,i_N\in [2]}\in \Sym(2,N)$.
 So $\|\cT\|_{\infty,2}=\|\cT'\|_{\infty,2}$.  Our assumption implies that there exists $\w'\in\rS^1$ such that $\|\cT'\|_{\infty,2}=
 |\an{\cT',\otimes_{j\in[N]}\w'_j}|$, where $\w_j'=\w'$ for $j\in [N]$.  Let $\w_j=\w'\oplus \0_{N-2}\in\rS^{N-1}$.  Then $\an{\cT,\otimes_{j\in[N]}\w_j}\otimes_{j\in[N]}\w_j$ is a symmetric best rank one approximation of $\cT$.  \qed

 \section{Best rank one approximations of $\cT\in\Sym(2,3)$}
 \begin{theo}\label{sym23case}  Let $\cT=[t_{i,j,k}]\in\Sym(2,3)$.  Then  each best rank one approximation of $\cT$ is symmetric,
 unless $\cT$ is a nonzero tensor proportional to the following one.
 \begin{equation}\label{sym23sten}
 t_{1,1,1}=\cos \theta, t_{1,1,2}=\sin \theta, t_{1,2,2}=-\cos\theta, t_{2,2,2}=-\sin\theta, \quad \theta\in [0,2\pi).
 \end{equation}
 For the above tensor there is a best rank approximation of the form $\uu\otimes\vv\otimes\w(\uu,\vv)$ where $\uu,\vv\in\rS^1$ are arbitrary
 and $\w(\uu,\vv)\in\rS^1$ is uniquely determined by $\uu,\vv$.  Furthermore, the above tensor has three symmetric best rank one approximations.
 \end{theo}
 \proof  Since $\cT=0$ has a unique best rank one approximation - $\cT$, we assume that $\cT\ne 0$.
 Suppose that $\cT$ has a best rank approximation $c\x\otimes\y\otimes \z, \x,\y,\z\in \rS^1$ which is not symmetric.
 (Note that $c\ne 0$.)
 By permuting $\x,\y,\z$ we can assume that $\y\ne \pm \z$.  Let $A(\x):=\cT\times \x\in\Sym(2,2)$.  (Because $\cT$ is symmetric, it is not important which index we contract.)  Hence $c\y\z\trans$ is best rank one approximation of the symmetric matrix $A(\x)$.  (Note that $A(\x)\ne 0$.) Proposition \ref{branapproxsymmat} yields that $\lambda_1(A(\x))+\lambda_2(A(\x))=0$, which is equivalent to $\tr (A(\x))=0$.
 Observe next that $A(\x)^2$ is a scalar matrix, i.e. $A(\x)^2=\lambda_1(A(\x))^2 I_2$.  Let $\uu\in\rS^1$ be arbitrary and $\vv:=\frac{1}{\|A(\x)\uu\|}A(\x)\uu(\in\rS^1)$.  Proposition \ref{branapproxsymmat} yields that $\|A(\x)\uu\|\uu\vv\trans$ is a best rank approximation of $A(\x)$.  Choose $\uu$ so that $\uu,\vv\not\in\{\x,-\x\}$.  Then $c'\x\vv\trans$ is a nonsymmetric best rank one approximation
 of $A(\uu)$.  The previous arguments show that $\tr (A(\uu))=0$.  Since $\x$ and $\uu$ are linearly independent, it follows that the two frontal
 section $A_k:=[t_{i,j,k}]_{i=j=1}^2\in \Sym(2,2)$ have trace zero.  Taking in account the $\cT\in\Sym(2,2)$ we deduce that $\cT$ is proportional
 to the tensor given by \eqref{sym23sten}.

 Assume that $\cT$ is of the form \eqref{sym23sten}.  Let
 \[A_1(\theta)=\left[\begin{array}{cc}\cos\theta&\sin\theta\\ \sin\theta&-\cos\theta\end{array}\right],  \quad
 A_2(\theta)=\left[\begin{array}{cc}\sin\theta&-\cos\theta\\ -\cos\theta&-\sin\theta\end{array}\right]\]
 be the two frontal sections of $\cT$.  Then
 \begin{equation}\label{Axform}
 A((\cos\phi,\sin\phi)\trans)=\cos\phi A_1(\theta)+\sin\phi A_2(\theta)=A_1(\theta-\phi).
 \end{equation}
 So $\lambda_1(A(\uu))=-\lambda_2(A(\uu))=1$ for every $\uu\in \rS^1$.  Hence any best rank one approximation of $A(\uu)$ is of the form
 $\vv\w(\uu,\vv)\trans$, where $\vv\in\rS^1$ and $\w(\uu,\vv):=A(\uu)\vv$.  This shows that any best rank approximation of $\cT$ is $\uu\otimes\vv\otimes\w(\uu,\vv)$ as claimed.

 It is left to show that $\cT$ has exactly $3$ different best rank one symmetric approximations.  In view of \eqref{Axform},
 by changing an orthonormal basis in $\R^2$ we may assume that $\theta=0$.  The condition that $\cT$ has a symmetric best rank approximation
 means that we need to choose $\x\in\rS^1$ such that $A(\x)=\pm \x$.  Note that if $A(\x)\x=\x$ then $A(-\x)(-\x)=-(-\x)$.
 Hence we need to find all $\x\in\rS^1$ such that $A(\x)\x=\x=(\cos\phi,\sin\phi)$.  This condition gives rise to the following three solutions
 \[(1,0)\trans, (-\frac{1}{2},\frac{\sqrt{3}}{2})\trans, (-\frac{1}{2},-\frac{\sqrt{3}}{2})\trans).\]
 \qed

 \section{Proofs of Theorems \ref{brank1symtheo} and \ref{maintheo}}
 \begin{lemma}\label{lem2dcase}
 Let $d\ge 2$ be an integer and $\cT=[t_{i_1,\ldots,i_d}]\in \Sym(2,d)\setminus\{0\}$.  Assume that $\cT$ has a nonsymmetric best rank one
 approximation.  Then for each $i_3,\ldots,i_d\in [2]$ the symmetric matrix $[t_{i,j,i_3,\ldots,i_d}]_{i=j=1}^2\in\Sym(2,2)$ has zero trace.
 \end{lemma}
 \proof  For $d=2$ the lemma follows from Proposition \ref{branapproxsymmat}.  For $d=3$ the lemma follows from Theorem \ref{sym23case}.
 We prove the lemma by induction on $d\ge 3$.  Suppose that the lemma holds for $d=N\ge 3$.  Assume that $d=N+1$.  Suppose that $\cA=a\otimes_{j=1}^d\x_j, \x_j\in\rS^1, j\in[d], a\ne 0$ is a nonsymmetric best rank one approximation of $\cT$.  Since $\cT$ is symmetric it follows that
 for each permutation $\sigma: [d]\to [d]$ the decomposable tensor $a\otimes_{j=1}^d \x_{\sigma(j)}$ is best rank one approximation of $\cT$ which is nonsymmetric.
 Hence, without a loss of generality we may assume that $\x_d\ne \pm\x_{d-1}$.  Fix the vectors $\x_1,\ldots,\x_{d-3}$ and consider
 $\cT(\x_1,\ldots,\x_{d-3}):=\cT\times \otimes_{j=1}^{d-3}\x_j\in\Sym(2,3)$.  Clearly, $a\x_{d-2}\otimes\x_{d-1}\otimes\x_d$ is a nonsymmetric
 best rank one approximation to $\cT(\x_1,\ldots,\x_{d-3})$.

 Theorem \ref{sym23case} yields that a best rank one approximation of $\cT(\x_1,\ldots,\x_{d-3})$ can be chosen of the form $a\otimes \w(\uu,\vv)\otimes \uu\otimes \vv$ where $\uu,\vv$ are arbitrary vectors in $\rS^1$ and $\w(\uu,\vv)\in\rS^1$.  Fix the vector $\vv=(v_1,v_2)\trans\in\rS^1$.
 Then $\cA:=a((\otimes_{j=1}^{d-3}\x_j)\otimes\w(\uu,\vv)\otimes\uu)$ is best rank one approximations of $\cT(\vv):=\cT\times \vv$.
 Observe that $\cA$ is not symmetric for $\uu\ne \pm \x_1$.
 Hence the induction hypothesis yield that the tensor $\cT(\vv)=[t_{i_1,\ldots,i_{d-1},1}v_1+t_{i_1,\ldots,i_{d-1},2}v_2]
 \in\Sym(2,d-1)$ satisfies that the assumption of the lemma.  I.e., for any $i_3,\ldots,i_{d-1}\in [2]$ the matrix $[t_{i,j,i_3,\ldots,i_{d-1},1}v_1+t_{i_1,\ldots,i_{d-1},2}v_2]_{i=j=1}^2$ has zero trace.  Let $\vv=(1,0)\trans,(0,1)\trans$ to deduce that
 the lemma holds for $\cT$.  \qed

 \textbf{Proof of Theorem \ref{brank1symtheo}}.   We first prove the theorem for $\cT=[t_{i_1,\ldots,i_d}]\in \Sym(2,d)\setminus\{0\}$.
 Suppose first that for some $i_3,\ldots,i_d\in[2]$ the $2\times 2$ symmetric matrix $[t_{i,j,i_3,\ldots,i_d}]_{i,j\in[2]}$ does not have trace $0$. Then Lemma \ref{lem2dcase} yields that each best rank one approximation tensor of $\cT$ is symmetric.

 Assume now that for each $i_3,\ldots,i_d\in[2]$ the $2\times 2$ symmetric matrix $[t_{i,j,i_3,\ldots,i_d}]_{i,j\in[2]}$ has trace $0$.
 Let $\cS=[s_{i_1,\ldots,i_d}]\in\Sym(2,d)$ be the following tensor.  $s_{1,\ldots,1}=1$ and all other entries of $\cS$ are zero.
 Then for any $\varepsilon\ne 0$ the tensor $\cT+\varepsilon\cS$ is symmetric, and the trace of the matrix $[t_{i,j,1,\ldots,1}+\varepsilon
 s_{i,j,1,\ldots,1}]_{i,j\in[2]}$ is $\varepsilon$.  Lemma \ref{lem2dcase} yields that $\cT+\varepsilon \cS$ has best rank one approximation
 $\cA(\varepsilon)\in \Sym(2,d)$.  Recall that
 \begin{equation}\label{Aepsprop}
 \|\cA(\varepsilon)\|^2=\|\cT+\varepsilon\cS\|_{2,\infty}^2=\an{\cT+\varepsilon \cS,\cA(\varepsilon)}.
 \end{equation}
 \eqref{pythid} yields that
 \[\|\cA(\varepsilon)\|\le \|\cT+\varepsilon \cA\|\le \|\cT\|+|\varepsilon|\|\cS\|=\|\cT\|+|\varepsilon|.\]
 Consider the bounded sequence $\cA_m:=\cA(\frac{1}{m})\in\Sym(2,d), m\in\N$.  There exists a subsequence $\cA_{m_i}, i\in\N$, such that
 $\lim_{i\to\infty}\cA_{m_i}=\cA\in\Sym(2,d)$, where $\cA$ is a decomposable tensor.  \eqref{Aepsprop} yields that $\cA$ is a best rank one
 approximation of $\cT$.  This completes the proof of the theorem for $\cT\in \Sym(2,d)$.

 Assume that $n>2$.  Use Lemma \ref{nto2} to conclude that each $\cT\in \Sym(n,d)\setminus\{0\}$ has a symmetric best rank one approximation.  \qed

 \textbf{Proof of Theorem \ref{maintheo}}.  Let $\cT\in \R^{n_1\times\ldots\times n_d}\setminus \{0\}$.  Assume that $[d]=\cup_{j=1}^k \alpha_j$ is the symmetric decomposition for $\cT$.  We prove the theorem by induction on $k$.  Assume that $k=1$, i.e. $\cT\in\Sym(n,d)$.  Then the theorem follows from Theorem \ref{brank1symtheo}.  Suppose that the theorem holds for $k\in [m]$, where $m \ge 1$ and assume that $k=m+1$.  Permute the factors in $\otimes_{j=1}^d\R^{n_i}$ to assume that $\alpha_1=\{1,\ldots,l\}, l<d$.  Suppose that $a\otimes_{j=1}^d \x_j, \x_j\in \rS^{n_j-1},j\in[d]$
 is a best rank one approximation of $\cT$.  Recall that $a\otimes_{j=l+1}^d \x_j$ is a best rank one approximation of
 $\cT(\x_1,\ldots,\x_l):= \cT\times \otimes_{j=1}^l \x_j$.  Furthermore \eqref{infshatnrmchar} yields that $\|\cT\|_{2,\infty}=\|\cT(\x_1,\ldots,\x_l)\|_{2,\infty}$.
 Observe that $\cT(\x_1,\ldots,\x_l)$ is symmetric with respect to  $\alpha_2,\ldots,\alpha_k$.  Hence the induction hypothesis yields that
 there exists a best rank one approximation $b\otimes_{j=l+1}^d \y_j$ of $\cT(\x_1,\ldots,\x_l)$, where $\y_j\in\rS^{n_j-1}$ for $j>l$, with the following properties.  $\y_p=\y_q$ for each $p,q\in\alpha_i, i>1$.  Let $\cT(\y_{l+1},\ldots,\y_k):=\cT\times\otimes_{j=l+1}^d\y_j$.
 Then $\cT(\y_{l+1},\ldots,\y_k)\in\Sym(n_1,l)$.  Theorem \ref{brank1symtheo} yields that $\cT(\y_{l+1},\ldots,\y_k)$ best rank one approximation
 of the form $c\otimes_{j=1}^l \y_j$ where $\y_1=\ldots=\y_l$.  Hence the rank one approximation $c\otimes_{j=1}^d\y_j$ is symmetric with respect
 $\alpha_1,\ldots,\alpha_k$.  \qed

 \section{Uniqueness of symmetric rank one approximation for generic symmetric tensors}
 For $\x=(x_1,\ldots,x_n)\trans \in \C^n$ denote $\|\x\|=\sqrt{\sum_{i=1}^n |\x_i|^2}$.  By $\otimes^d\C^n$ we denote the tensor products of $d$ copies of $\C^n$.  Then $\otimes^d\x\in\otimes^d \C^n$ is a decomposable tensor $\x\otimes\ldots\otimes\x$.
 Denote by $\Sym(n,d,\C)\subset \otimes^d\C^n$ the space of all symmetric tensors $\cT=[t_{i_1,\ldots,i_d}]_{i_1=\ldots=i_d=1}^n$ with
 complex entries.  It is well known that $\Sym(n,d)$ is isomorphic to $\C^{n+d-1\choose d}$.

 \begin{theo}\label{numbcritsymten}  Let $d\ge 3, n\ge 2$ be integers.  Then there exists an algebraic variety $\Sigma(n,d)\subset\Sym(n,d,\C)$ such that for each $\cT\in \Sym(d,n,\C)\setminus \Sigma(n,d)$ the symmetric eigensystem
 \begin{equation}\label{symeigensys}
 \cT\times \otimes^{d-1}\x=\x, \quad \x\ne \0
 \end{equation}
 have at most  $(d-1)^n-1$ distinct solutions, which are invariant under the multiplication by $d-2$ roots of unity.

 Assume furthermore that $\cT\in\Sym(d,n)\setminus\Sigma(n,d)$.   If $d\ge 3$ is odd,  the number of real distinct solutions of the above
 system is at most $\frac{(d-1)^n-1}{d-2}$.  Furthermore any solution of
 \begin{equation}\label{symeigensys1}
 \cT\times \otimes^{d-1}\x=-\x, \quad \x\ne \0
 \end{equation}
 is the negative of the real solution of \eqref{symeigensys}.

 If $d$ is even then the systems \eqref{symeigensys}--\eqref{symeigensys1} all together have at most $\frac{2((d-1)^n-1)}{d-2}$ real solutions which are invariant by multiplication by $-1$.
 \end{theo}
 \proof  Let $\cT\in \Sym(n,d,\C)$. Consider the system
 \begin{equation}\label{homeigsys}
 \cT\times \otimes^{d-1}\x=\0,\quad \x\in\C^n.
 \end{equation}
 Let $\uu_1,\ldots,\uu_n\in \C^n$ be linearly independent.  Let $\cT_0=\sum_{i=1}^n \otimes^d \uu_i$.  Then the above system is equivalent to
 $\sum_{i=1}^d (\uu_i\trans \x)^{d-1}\uu_i=\0$.  Since $\uu_1,\ldots,\uu_n$ are linearly independent it follows that
 $(\uu_i\trans \x)^{d-1}=0$ for $i\in[n]$.  Hence $(\uu_i\trans \x)=0$ for $i\in[n]$.  As $\uu_1,\ldots,\uu_n$ independent it follows that $\x=\0$.
 That is
 \[\min_{\|\x\|=1,\x\in\C^n} \|\cT_0\times \otimes^{d-1}\x\|>0.\]
 The above inequality holds for $\cT\in\Sym(n,d,\C)$ in some neighborhood of $\cT_0$.
 Hence there exists a strict algebraic variety $\Sigma(n,d)\subset\Sym(n,d,\C)$ such that for each $\cT\in \Sym(n,d,\C)\setminus \Sigma(n,d)$
 the system \eqref{homeigsys} has a unique solution $\x=0$.

 Assume that $\cT\in \Sym(n,d,\C)\setminus \Sigma(n,d)$.  Consider the polynomial system
 \begin{equation}\label{symeigensys0}
 \cT\times \otimes^{d-1}\x-\x=\0, \quad \x\in\C^n.
 \end{equation}
 Its principal part is the homogeneous system \eqref{homeigsys}, where each homogeneous polynomial
 is of degree $d-1$.  Hence \eqref{symeigensys0} has at most $(d-1)^n$ distinct solutions.  (This is the precise version of Bezout's theorem.  See for example \cite{Fr77}.)
 Note that $\x=\0$ is a solution.  Hence \eqref{symeigensys} has at most $(d-1)^n-1$ distinct solutions.
 Clearly, if $\x\ne \0$ is a solution of \eqref{symeigensys} then $\zeta\x$ is also a solution of \eqref{symeigensys} if $\zeta^{d-2}=1$.
 As $\x\ne 0$, note that $\zeta\x\ne \eta\x$ if $\zeta,\eta$ are two distinct $d-2$ roots of $1$.  That establishes the first part of the
 of the theorem.

 Assume now that $\cT\in \Sym(n,d)\setminus\Sigma(d,n)$.  We know that \eqref{symeigensys} has at most $(d-1)^n-1$ complex distinct solutions.
 Assume that $\x\ne \0$ is a real solution.  Then $\zeta\x$ is also a solution for $\zeta^{d-2}=1$.
 Suppose that $d\ge 3$ is odd.  Then $\zeta\x$ is real when $\zeta^{d-2}=1$ if and only if $\zeta=1$.
 Hence each real solution of \eqref{symeigensys} gives rise to another $d-3$ distinct nonreal solutions of \eqref{symeigensys}.
 Hence the number of real solutions of \eqref{symeigensys} is at most  $\frac{(d-1)^n-1}{d-2}$.   Note that for any real solution $\x$
 of \eqref{symeigensys} $-\x$ is a real solution of \eqref{symeigensys1} and vice versa.  This proves our theorem for $d\ge 3$ odd.

 Assume now that $d\ge 3$ is even.  Then for any real solution $\x$ of \eqref{symeigensys} or \eqref{symeigensys1} we get another real solution $-\x$, and $d-4$ complex solutions of \eqref{symeigensys} or \eqref{symeigensys1}, respectively, of the form $\zeta\x$, where $\zeta^{d-2}=1,\zeta\ne \pm 1$.  Assume now that $\x$ is a real solution \eqref{symeigensys1}.
 Then $\eta\x, \eta^{d-2}=-1$ gives rise to $d-2$ distinct nonreal solutions of \eqref{symeigensys}.  Hence the total number of real solutions of the systems
 \eqref{symeigensys}--\eqref{symeigensys1} is $\frac{2((d-1)^n-1)}{d-2}$.  \qed
 \begin{theo}\label{simcritsymtens}  There exists a subvariety $\Sigma_1\subset\Sym(n,d,\C)$, where $\Sigma_1(n,d)\supseteq \Sigma(n,d)$, with the following
 properties.  Let $\cT\in\Sym(n,d)\setminus\Sigma_1(d,n)$ and consider nonzero critical values $\an{\cT,\otimes^d\x}$ on $\rS^{n-1}$ and the corresponding
 critical points.  Then
 \begin{enumerate}
 \item\label{simcritsymtens1}
  For an odd $d\ge 3$ there are $N(\cT)\le \frac{2((d-1)^n-1)}{d-2}$ nonzero distinct critical values, where each critical value $\lambda$ has a unique corresponding critical point $\x\in\rS^{n-1}$.  Furthermore, $-\lambda$ is also a critical value with the corresponding unique critical point $-\x$.
 \item\label{simcritsymtens2}
  For an even $d\ge 3$ there are $N(\cT)\le \frac{2((d-1)^n-1)}{d-2}$ nonzero distinct critical values, where each critical value $\lambda$ has exactly two critical points $\pm\x\in\rS^{n-1}$.  Furthermore, the absolute values of two distinct critical points are distinct.
 \end{enumerate}
 \end{theo}
 \proof  Let $\cT\in \Sym(n,d,\C)\setminus\Sigma(n,d)$.  We claim that there exists a subvariety $\Sigma_0(d,n)\subset
 \Sym(n,d,\C)$ such that for each $\cT\not\in\Sigma_0(n,d)$ the following conditions hold.
 \begin{enumerate}
 \item\label{doddsig1con}
 Assume that $d\ge 3$ is odd.  Then the system \eqref{symeigensys} has $(d-1)^n-1$ solutions, and for two different solutions $\x,\y$ the inequality $\x\trans  \x\ne \y\trans \y$ holds.
 \item\label{devensig1con}
 Assume that $d\ge 3$ is even.  Then the systems \eqref{symeigensys} and \eqref{symeigensys1} each has $(d-1)^n-1$ solutions, and for two different solutions $\x,\y$ of either \eqref{symeigensys} or \eqref{symeigensys1},  where $\y\ne -\x$,
 the inequality $\x\trans  \x\ne \y\trans\y$ holds.
 \end{enumerate}

 We consider a special tensor $\cT=\sum_{i=1}^nt_i^{-(d-2)}\otimes^d \e_i$, where $\e_i=(\delta_{1i},\ldots,\delta_{ni})\trans, t_i\in\C\setminus\{0\}$ for $i\in[n]$.
 For this tensor we can explicitly calculate all solutions of \eqref{symeigensys}.
 Namely, it is of the form $(x_1,\ldots,x_n)\trans$ where each $\frac{x_i}{t_i}$ satisfies the equation $x(x^{d-2}-1)=0$.
 (We need to exclude from this list the trivial solution $\x=\0$.)
 This means that \eqref{symeigensys} has  $(d-1)^n-1$ solutions.  Moreover, the solutions of \eqref{symeigensys1}
 are of the form $\theta\x$, where $\x$ is a solution of  \eqref{symeigensys} and $\theta$ is a fixed solution of $\theta^{d-2}=-1$.
 To be explicit, let $m\ge 0$ be the integer such that $\frac{d-2}{2^m}$ is an odd integer.  Then $\theta:=e^{\frac{\pi\sqrt{-1}}{2^{m}}}$.

 Let $\F=\Q[\xi]$ be a finite extension field of the rational numbers $\Q$, where $\xi$ is a primitive root of $\xi^{d-2}=1$.
 (Each element of $\F$ is a polynomial of degree $d-3$ at most with rational coefficients.)
 Assume that $t_1^2,\ldots,t_n^2$ are linearly independent over $\Q[\xi]$.

 We first consider the simple case: $d\ge 3$ is odd.
 We claim that if $\x=(x_1,\ldots,x_n)\trans$ and $\y=(y_1,\ldots,y_n)\trans$ are different nonzero solution of \eqref{symeigensys}  then $\x\trans\x\ne \y\trans\y$.
 Indeed $\x=(\zeta_1 t_1,\ldots,\zeta_n t_n)\trans, \y=(\eta_1 t_1,\ldots,\eta_n t_n)\trans$.
 The assumption that $\x$ and $\y$ satisfy \eqref{symeigensys} imply that each $\zeta_i$ and $\eta_i$
 satisfy the equation $s(s^{d-2}-1)=0$.  Observe next that
 \begin{equation}\label{xyeq}
 \x\trans\x-\y\trans\y=\sum_{i=1}^n (\zeta_i^2-\eta_i^2)t_i^2.
 \end{equation}

 Suppose that $\x\trans\x-\y\trans\y=\0$.  As $t_1^2,\ldots,t_n^2$ are linearly independent over $\F$ it follows that
 $\zeta_i^2=\eta_i^2$ for $i\in[d]$.  So $\zeta_i=\pm \eta_i$ for $i\in [d]$.  Clearly $\zeta_i=0\iff \eta_i=0$.
 Assume that $\zeta_i\ne 0$.  So $\zeta_i^{d-2}=1\Rightarrow \eta_i^{d-2}=1$.  Hence $\zeta_i=\eta_i$ for $i\in [d]$.
 So $\x=\y$ contrary to our assumption.  Hence, there exists a subvariety $\Sigma_0(n,d)\subset \Sym(n,d,\C)$ such that
 for $\cT\in \Sym(n,d,\C)\setminus(\Sigma(n,d)\cup \Sigma_0(n,d))$ the condition \ref{doddsig1con} hold.

 Assume now that $d\ge 3$ is even.  Let $\x=(\zeta_1 t_1,\ldots,\zeta_n t_n)\trans$ and $\y=(\eta_1 t_1,\ldots,\eta_n t_n)\trans$
 satisfy either \eqref{symeigensys} or \eqref{symeigensys1}.  So $\zeta_i((\phi\zeta_i)^{d-2})=\eta_i((\psi\eta_i)^{d-2}-1)=0$ for $i\in [d]$.
 Here $\phi,\psi\in\{1,\theta\}$.  Assume that $\x\trans\x-\y\y\trans=0$.  So $\zeta_i=\pm \eta_i$ for $i\in [d]$.  Hence $\phi=\theta$.
 That is either $\x$ and $\y$ satisfy \eqref{symeigensys} or $\x$ and $\y$ satisfy \eqref{symeigensys1}.
 However, it is possible that $\x\trans\x-\y\trans\y=0$ and $\x\ne \pm \y$.  We now find a tensor $\cT'\in \Sym(n,d,\C)$ in a small
 neighborhood of $\cT$ such that for $\x',\y'$, where $\x'\ne \pm\y'$, satisfying either the system \eqref{symeigensys} or \eqref{symeigensys1}
 for $\cT'$, one has the inequality $(\x')\trans\x'-(\y') \trans\y'\ne 0$.  Since any solution $\x$ of \eqref{symeigensys} or \eqref{symeigensys1}
 is a simple solution, (as we have the maximal number of distinct solutions), we can use an implicit function theorem to find the unique solutions $\x'$ of \eqref{symeigensys} or \eqref{symeigensys1} in the neighborhood of the solution $\x$ respectively.  Since for $\x\trans\x-\y\trans\y=0$
 may hold if only $\x,\y$ satisfy both either \eqref{symeigensys} or \eqref{symeigensys1}  it is enough to show that $(\x')\trans\x'-(\y') \trans\y'\ne 0$ where $\x',\y'$ satisfy \eqref{symeigensys}.

 Let $\cT(\varepsilon):=\cT+\varepsilon\cS$ for a fixed $\cS\in \Sym(n,d,\C)$.  Consider the system
 \begin{equation}\label{epseigsys}
 \cT(\varepsilon)\times\otimes^{d-1}\x(\varepsilon)=\x(\varepsilon).
 \end{equation}
 The implicit function theorem yields that $\x(\varepsilon)$ can be expanded in power series of $\varepsilon$.
 So
 \begin{equation}\label{powexpan}
 \x(\varepsilon)=\x+\varepsilon \x_1+O(\varepsilon^2).
 \end{equation}
 Denote by $\cT(\x):=\cT\times\otimes^{d-2}\x\in \Sym(n,2,\C)$.  The system \eqref{symeigensys} is equivalent to $\cT(\x)\x=\x$ and
 $\x\trans \cT(\x)=\x\trans$.
 Then the $\varepsilon$ term in \eqref{epseigsys} is
 \begin{equation}\label{ftermepeig}
 (d-1)\cT(\x)\x_1+\cS\times \otimes^{d-1}\x=\x_1.
 \end{equation}

 Multiply the above equality by $\x\trans$ to deduce that
 \begin{equation}\label{xx1teq}
 \x\trans\x_1=\frac{\an{\cS,\otimes^d\x}}{2-d}.
 \end{equation}
 Observe finally that
 \begin{equation}\label{xepyepdif}
 \x(\varepsilon)\trans\x(\varepsilon)-\y(\varepsilon)\trans\y(\varepsilon)=\x\trans\x-\y\trans\y +\frac{2\varepsilon}{2-d}\an{\cS,\otimes^d\x-
 \otimes^d\y}+O(\varepsilon^2).
 \end{equation}

 Note that $\otimes^d\x- \otimes^d\y$ is zero tensor if and only if $\x=\gamma\y$ where $\gamma^d=1$.
 Hence we can choose $\cS\in \Sym(n,d,\C)$ lying outside a finite number of subspaces of codimension one, such that
 the coefficient of $\varepsilon$ is different from zero for each pair of solutions
 $\x,\y$ of \eqref{symeigensys} such that $\x\ne \pm\y$ and $\x\trans\x=\y\trans\y$.
 Hence, there exists a subvariety $\Sigma_0(n,d)\subset \Sym(n,d,\C)$ such that
 for $\cT\in \Sym(n,d,\C)\setminus(\Sigma(n,d)\cup \Sigma_0(n,d))$ the condition \ref{devensig1con} hold.

 Let $\Sigma_1(n,d)=\Sigma(n,d)\cup\Sigma_0(n,d)$ and
 $\cT\in \Sym(n,d)\setminus \Sigma_1(n,d)$.  Assume that $\lambda\ne 0$ is a critical value of $\an{\cT,\otimes^d \z}$
 on $\rS^{n-1}$ with a corresponding critical point $\y\in\rS^{n-1}$.  Hence $\cT\times\otimes^{d-1}\y=\lambda\y$.

 Assume first that $d\ge 3$ is odd.
 Clearly, $-\y$ is a critical point corresponding to the critical value $-\lambda$.
 Without loss of generality we may assume that $\lambda>0$.  Then $\x:=\lambda^{-\frac{1}{d-2}}\y$  is a real solution of \eqref{symeigensys}.
 Vice versa, any real solution of $\x$ of \eqref{symeigensys} gives rise to a critical point $\y=\frac{1}{\|\x\|}\x$ with the critical value
 $\|\x\|^{-(d-2)}$.  Recall that for $\cT\in\Sym(n,d)\setminus\Sigma_1(n,d)$ \eqref{symeigensys} has exactly $(d-1)^n-1$ different complex solutions,
 with corresponding $(d-1)^n-1$ different values of $\x\trans\x$.  Hence \eqref{symeigensys} has at most $\frac{(d-1)^n-1}{d-2}$ real solutions
 and the length of all these solutions are distinct.  Therefore, the number of positive nonzero critical points of $\an{\cT,\otimes^d\y}$ on $\rS^{n-1}$ is at most $\frac{(d-1)^n-1}{d-2}$, and to each critical positive value corresponds a unique critical point.
 These arguments prove the theorem for $d$ odd.

 Assume now that $d$ is even.  The $-\y$ is a critical point of the critical
 value $\lambda$.  Let $\x_+:=\lambda^{-\frac{1}{d-2}}\y$ if $\lambda>0$ and $\x_-:=(-\lambda)^{-\frac{1}{d-2}}\y$ if $\lambda<0$.
 Then $\x_+$ satisfies \eqref{symeigensys} and $\x_-$ satisfies \eqref{symeigensys1}.
 Theorem \ref{numbcritsymten} yields that the number of critical nonzero values of $\an{\cT,\otimes^d \z}$ is at most $\frac{2((d-1)^n-1)}{d-2}$.
 It is left to show that the absolute values of nonzero critical values are distinct.
 Let $\x,\y$ be two solutions of either \eqref{symeigensys} or \eqref{symeigensys1}, and assume that $\x\ne \pm \y$.
 Then $\uu:=\frac{1}{\|\x\|}\x,\vv:=\frac{1}{\|\y\|}\y$ are critical points of $\an{\cT,\otimes^d \z}$ on $\rS^{n-1}$, corresponding
 to critical values $\lambda,\mu$ respectively.  Clearly $|\lambda|=\|\x\|^{-(d-2)}, |\lambda|=\|\y\|^{-(d-2)}$.
 Since $\cT\not\in \Sigma_1(n,d)$ we deduce that $\x\trans\x\ne \y\trans\y$.  So $|\lambda|\ne |\mu|$.  \qed

 \begin{corol}\label{uniquerank1symap}  Let $d\ge 3$ be an integer and assume that $\cT\in\Sym(n,d)\setminus\Sigma_1(n,d)$.  Then $\cT$ has a unique symmetric best rank one approximation.
 \end{corol}
 \proof  Best rank one approximation corresponds to the critical point $\y$ of $\an{\cT,\otimes^d\z}$, corresponding to the critical
 value $\lambda$ where $\lambda^2$ is maximal.  The corresponding best rank one symmetric approximation is $\cA:=\an{\cT,\otimes^d\y}\otimes^d\y$.
 So $-\y$ gives the same $\cA$.  For $d$ odd all positive critical values are distinct.  Hence $\cA$ is unique.
 For $d$ even, each nonzero critical value has two critical points $\pm \y$.  The absolute values of the critical values are distinct.
 Hence $\cA$ is unique.  \qed

 \bibliographystyle{plain}

\begin{thebibliography}{MMM}

 \bibitem{CHLZ12} B. Chen, S. He, Z. Li, and S. Zhang, Maximum block improvement and polynomial optimization,
 SIAM J. OPTIM. 22 (2012),  87–-107.

 \bibitem{DF93} A. Defant and K. Floret, \emph{Tensor Norms and Operator Ideals}, North-Holland, Amsterdam, 1993.

 \bibitem{Fr77} S. Friedland, Inverse eigenvalue problems, {\it Linear Algebra Appl.} 17
 (1977), 15-51.

 \bibitem{Fr82} S. Friedland,  Variation of tensor powers and spectra, \emph{Linear
 Multilin. Algebra} 12 (1982), 81-98.

 \bibitem{FO12}  S. Friedland and G. Ottaviani, The number of singular vector
 tuples and uniqueness of best rank one approximation of tensors,
 arXiv:1210.8316.

 \bibitem{FRS}  S. Friedland, J. Robbin and J. Sylvester,  On the crossing rule,
 \emph{Comm. Pure Appl. Math.} 37 (1984), 19-37.


 \bibitem{GolV96} G.H. Golub and  C.F. Van Loan, {\it Matrix
 Computation}, John Hopkins Univ. Press, 3rd Ed., 1996.


 \bibitem{Kol} T. Kolda,  On the best rank-$k$ approximation of a symmetric tensor, Householder 2011.


 \bibitem{Lim05} L.-H. Lim, Singular values and eigenvalues of tensors:
 a variational approach, \emph{Proc. IEEE International Workshop on
 Computational Advances in Multi-Sensor Adaptive
 Processing} (CAMSAP '05), 1 (2005), 129-132.













 \end{thebibliography}
 
 \end{document}